\documentclass[a4paper,12pt]{amsart}
\usepackage{amssymb}
\usepackage{ifthen}
\usepackage{cite}
 \usepackage[dvips]{graphicx}
\nonstopmode\numberwithin{equation}{section}
\usepackage{amssymb}
\usepackage{ifthen}
\usepackage{graphicx}
\usepackage{amsmath}
\usepackage[T1]{fontenc} 
\usepackage[utf8]{inputenc}
\usepackage[usenames,dvipsnames]{color}
\usepackage{color}
\usepackage[english]{babel}
\usepackage{fancyhdr}
\usepackage{fancybox}
\usepackage{tikz}

\usepackage{color,xcolor}

\newtheorem*{theoA}{Theorem A}
\newtheorem*{theoB}{Theorem B}
\newtheorem*{theoC}{Theorem C}

\nonstopmode \numberwithin{equation}{section}
\theoremstyle{plain}

\newtheorem{conj}{Conjecture}

\theoremstyle{definition}
\newtheorem{defi}{Definition}[section]

\newtheorem{thm}{Theorem}[section]
\newtheorem{prob}{Problem}[section]
\newtheorem{cor}{Corollary}[section]
\newtheorem{ques}{Question}[section]
\newtheorem{prop}{Proposition}[section]
\newtheorem{rem}{Remark}[section]
\newtheorem{lem}{Lemma}[section]

\newtheorem{Open Ques}{Open question}[section]

\newcounter{minutes}
\divide\time by 60
\newcounter{hours}
\multiply\time by 60
\addtocounter{minutes}{-\time}

\newcounter {own}
\def\theown {\thesection       .\arabic{own}}

\newenvironment{pf}[1][]{%
 \vskip 3mm
 \noindent
 \ifthenelse{\equal{#1}{}}%
  {{\slshape Proof. }}%
  {{\slshape #1.} }%
 }%
{\qed\bigskip}

\newcounter{alphabet}

\def\be{\begin{equation}}
\def\ee{\end{equation}}

\newcommand{\bee}{\begin{enumerate}}
\newcommand{\eee}{\end{enumerate}}

\newcommand{\blem}{\begin{lem}}
\newcommand{\elem}{\end{lem}}
\newcommand{\bthm}{\begin{thm}}
\newcommand{\ethm}{\end{thm}}
\newcommand{\bcor}{\begin{cor}}
\newcommand{\ecor}{\end{cor}}
\newcommand{\beg}{\begin{examp}}
\newcommand{\eeg}{\end{examp}}
\newcommand{\begs}{\begin{examples}}
\newcommand{\eegs}{\end{examples}}
\usepackage{mathrsfs}
\newcommand{\bdefn}{\begin{defn}}
\newcommand{\edefn}{\end{defn}}

\newcommand{\bprob}{\begin{prob}}
\newcommand{\eprob}{\end{prob}}
\newcommand{\bei}{\begin{itemize}}
\newcommand{\eei}{\end{itemize}}

\newcommand{\bcon}{\begin{conj}}
\newcommand{\econ}{\end{conj}}
\newcommand{\bcons}{\begin{conjs}}
\newcommand{\econs}{\end{conjs}}
\newcommand{\bprop}{\begin{prop}}
\newcommand{\eprop}{\end{prop}}
\newcommand{\br}{\begin{rem}}
\newcommand{\er}{\end{rem}}
\newcommand{\brs}{\begin{rems}}
\newcommand{\ers}{\end{rems}}
\newcommand{\bo}{\begin{obser}}
\newcommand{\eo}{\end{obser}}
\newcommand{\bos}{\begin{obsers}}
\newcommand{\eos}{\end{obsers}}
\newcommand{\bpf}{\begin{pf}}
\newcommand{\epf}{\end{pf}}
\newcommand{\ba}{\begin{array}}
\newcommand{\ea}{\end{array}}
\newcommand{\beq}{\begin{eqnarray}}
\newcommand{\beqq}{\begin{eqnarray*}}
\newcommand{\eeq}{\end{eqnarray}}
\newcommand{\eeqq}{\end{eqnarray*}}

\begin{document}

\title{Bohr phenomena for slice regular functions over Quaternions}

\author{Sabir Ahammed}
\address{Sabir Ahammed, Department of Mathematics, Jadavpur University, Kolkata-700032, West Bengal, India.}
\email{sabira.math.rs@jadavpuruniversity.in}

\author{Molla Basir Ahamed}
\address{Molla Basir Ahamed, Department of Mathematics, Jadavpur University, Kolkata-700032, West Bengal, India.}
\email{mbahamed.math@jadavpuruniversity.in}

\author{Ming-Sheng Liu$^*$}
\address{Ming-Sheng Liu, School of Mathematical Sciences, South China Normal University, Guangzhou, Guangdong 510631, China.}
\email{liumsh65@163.com}

\subjclass[{AMS} Subject Classification:]{Primary 30C35, 30C50, 30G35, 30C99,  17A35.}

\def\thefootnote{}
\footnotetext{ {\tiny File:~\jobname.tex,
printed: \number\year-\number\month-\number\day,
          \thehours.\ifnum\theminutes<10{0}\fi\theminutes }
} \makeatletter\def\thefootnote{\@arabic\c@footnote}\makeatother

\maketitle
\pagestyle{myheadings}
\markboth{S. Ahammed, M. B. Ahamed, and M.-S. Liu}{Bohr phenomena for slice regular functions}
\begin{abstract}
Slice regular functions are a generalization of holomorphic functions to the setting of quaternions (and more generally, Clifford algebras). In this paper, we first establish the Bohr inequality for slice starlike functions and slice close-to-convex functions over quaternions $\mathbb{H}$. Next, we present a generalization of the Bohr inequality, and improved versions of the Bohr inequality for slice regular functions on the open unit ball $\mathbb{B}$ of $\mathbb{H}$. Finally, we provide a refined version of the Bohr inequality for slice regular functions $f$ on $\mathbb{B}$ such that $  {\rm Re}(f(q)) \leq 1 $ for all $q \in \mathbb{B}$. All the results are demonstrated to be sharp.
\end{abstract}
\section{\bf Introduction}
In the early twentieth century, the study of Dirichlet series of the form $\sum_{k=0}^{\infty}a_kk^{-s}$, where $a_k$ and $s$ are the complex numbers, became a highly attractive research topic. Harald Bohr initiated investigations into their convergence, focusing primarily on determining the largest region in which a given Dirichlet series converges absolutely. The absolute convergence problem of Dirichlet series (see \cite[p. xxiv]{Defant-2019}) aims to determine the maximum possible width of the strip where the series converges uniformly but not absolutely. Equivalently, this problem can be reformulated as finding the supremum of the abscissa of absolute convergence for a Dirichlet series that converges on $\{z\in\mathbb{C} : {\rm Re}(z)>0\}$ with a bounded limit function.\vspace{2mm}

This study is related to the classical result known as Bohr's theorem, which applies to the class $\mathcal{H}(\mathbb{D})$ of analytic functions $f(z) = \sum_{k=0}^{\infty} a_k z^k$ in the open unit disk $\mathbb{D}:=\{ z \in \mathbb{C} : |z| < 1 \}$, where $|f(z)| \leq 1$ for all $z \in \mathbb{D}.$  The result Harald Bohr first presented in 1914 has been improved by many renowned mathematicians. It has since become a foundation for research in numerous different function spaces. The Bohr theorem, in its ultimate form, is as follows.
\begin{theoA}{\rm \cite{Bohr-1914}} If $f(z)=\sum_{k=0}^{\infty}a_kz^k\in \mathcal{H}(\mathbb{D}),$  then
	\begin{align}\label{Eq-1.11}
		\sum_{k=0}^{\infty}|a_k||z|^k\leq 1 \;\;\;\mbox{for}\;\;\; |z|\leq \dfrac{1}{3}.
	\end{align}
	The constant $1/3$ is best possible.
\end{theoA}
The inequality $\eqref{Eq-1.11}$ fails when $ |z|>{1}/{3} $ in the sense that there are functions in $ \mathcal{H}\left(\mathbb{D}\right) $ (e.g., $f_a(z)=(a-z)/(1-az)$, where $a\in\mathbb{D}$) for which the inequality is reversed when $ |z|>{1}/{3} $. H. Bohr initially showed in \cite{Bohr-1914} the inequality \eqref{Eq-1.11} holds only for $|z|\leq1/6$, which was later improved independently by M. Riesz, I. Schur, F. Wiener and some others. The sharp constant $1/3$  and the inequality \eqref{Eq-1.11} in Theorem A are called respectively, the Bohr radius and the  classical Bohr inequality for the family $ \mathcal{H}\left(\mathbb{D}\right) $.  Several other proofs of this interesting inequality were given in different articles  (see \cite{Sidon-1927,Paulsen-Popescu-Singh-PLMS-2002,Tomic-1962}). Paulsen \emph{et al.} \cite{Paulsen-Popescu-Singh-PLMS-2002} were the first to demonstrate that the Bohr radius depends on the power of the term $|f(0)|$. They proved that the Bohr radius is $1/2$ when $|f(0)|$ is replaced by $|f(0)|^m$, where $m=2$. Blasco \cite{Blasco-OTAA-2010} extended the range of $m \in [1,2]$ and established a refinement of the classical Bohr inequality on $\mathbb{D}$.\vspace{1.2mm}

In $1995,$ the Bohr inequality was applied by Dixon \cite{Dixon-BLMS-1995} to the long-standing problem of characterizing Banach algebras that satisfying the von Neumann inequality and then this application garnered considerable attention from various researchers, leading to investigations of the phenomenon in a variety of function spaces. However, it can be noted that not every class of functions has Bohr phenomenon, for example, B\'en\'eteau \emph{et al.} \cite{BDK2004} showed that there is no Bohr phenomenon in the Hardy space $ H^p(\mathbb{D},X), $ where $p\in [1,\infty).$ Moreover, $\mathcal{H}\left(\mathbb{D}\right)$ is not the only class of analytic functions where the Bohr radii are studied; many other classes of functions and some integral operators also have their Bohr radii examined. The study of the Bohr phenomenon for holomorphic functions in several complex variables and in Banach spaces has becomes very active (see \cite{Aizenberg-PAMS-2000, Galicer-Mansilla-Muro-TAMS-2020, Hamada-IJM-2009, Hamada-RM-2024, Kumar-PAMS-2022,Kumar-Manna-JMAA-2023, Lata-Singh-PAMS-2022, Lin-Liu-Ponnusamy-Acta-2023, Liu-Ponnusamy-PAMS-2021}); in the non-commutative setting ( see \cite{Paulsen-Singh-BLMS-2006,Popescu-IEOT-2019}),  \cite{Paulsen-Popescu-Singh-PLMS-2002} for the operator-valued generalizations of the Bohr theorem in the single-variable case, see \cite{ Popescu-TAMS-2007} and references therein. In \cite{Kayumov-Khammatova-Ponnusamy-CRMA-2020} the authors have studied the Bohr radius for Ces\'aro operator. In \cite{Allu-Halder-PEMS-2023}, the authors have studied the Bohr inequality for operator-valued holomorphic functions defined in $\mathbb{D}.$ In \cite{Hamada-Math.Nachr-2021}, the author have the Bohr inequality for holomorphic and pluriharmonic mappings with values in complex Hilbert spaces. In \cite{Blasco-CM-2017}, the author studied the p-Bohr radius of a Banach space.\vspace{1.2mm}

An interesting application of the classical Bohr inequality for the class  $\mathcal{H}(\mathbb{D})$ was explored by Ponnusamy \emph{et al.} \cite[Theorem 2]{Ponnusamy-Vijayakumar-Wirths-RM-2020}, who investigated a modified version of the classical Bohr theorem, known as a refinement of the classical Bohr inequality. They present the following sharp inequality as follows:
\begin{align}\label{BS-Eq-1.2}
	\sum_{k=0}^{\infty}|a_kz^k|+\left(\frac{1}{1+|a_0|}+\frac{|z|}{1-|z|}\right)\sum_{k=1}^{\infty}|a_kz^{k}|^2\leq 1 \quad \mbox{for} \quad |z| \leq \frac{1}{2+|a_0|}
\end{align}
and the constant ${1}/{(2+|a_0|)}$ is sharp. \vspace{1.2mm}

As a generalization of holomorphic functions of one complex variable, the theory of slice regular functions of one quaternionic variable was initiated by Gentili and Struppa \cite{Gentili-Struppa-AVM-2007} and further developed by Clifford algebras \cite{Colombo-Sabadini-Struppa-IJM-2009} and octonions \cite{Gentili-Struppa-RMJM-2010}. Liu et al. \cite{Kou-Liu-Zou-AACA-2020, Liu-Xiang-Yang-AACA-2017} studied the Plancherel theorems of quaternion Hilbert transforms and the zeroes of Clifford algebra-valued polynomials. Based on the concept of stem functions, these three classes of functions were eventually unified and generalized into real alternative algebras \cite{Ghiloni-Perotti-AdvM-2011}. Following the historical path, the Bohr inequality has recently been generalized for slice regular functions in \cite{Xu-PRSE-2021} over the largest alternative division algebras of octonions.  Rocchetta \emph{et al.}~\cite{Rocchetta-Gentili-Sarfatti-MN-2012} established Bohr’s theorem for a class of regular functions over the quaternions, making extensive use of regular quaternionic rational transformations as developed in \cite{Bisi-Gentili-IUMJ-2009,Stoppato-AG}. This approach is made possible by the fact that quaternions form a skew field that is still associative. When restricted to each slice, the function \( f \) described in Theorem~B can be interpreted, following~\cite[Lemma 3.2]{Xu-PRSE-2021}, as a vector-valued holomorphic function from the unit disk \( \mathbb{D} \) into the unit ball of \( \mathbb{C}^4 \), equipped with the standard Euclidean norm. Blasco~\cite[Theorem 1.2]{Blasco-OTAA-2010} introduced the Bohr radius for holomorphic functions mapping \( \mathbb{D} \) into the unit ball of \( \mathbb{C}^n \) (for \( n \geq 2 \)) and showed that this radius is zero. From this perspective, the theory of slice regular functions differs from that of vector-valued holomorphic functions.
\vspace{2mm}

Slice regular functions provide a powerful and natural extension of complex holomorphic functions to quaternions. They preserve many classical results while allowing for new applications in physics, engineering, and functional analysis. Many classical theorems in complex analysis, such as the Identity Theorem, Maximum Modulus Principle, and Liouville’s Theorem, have analogues in the slice regular setting. Their structure is significantly influenced by the non-commutativity of quaternions, making their study both challenging and interesting.
\vspace{2mm}

Let $\mathcal{SB}(\mathbb{B})$ denote the class of slice regular functions  $f(q)=\sum_{k=0}^{\infty}q^kp_k$ on the open unit ball $\mathbb{B}$  of quaternions $\mathbb{H}$,  continuous on the closure $\overline{\mathbb{B}} ,$ such that $|f(q)|<1$ for all $|q|\leq 1.$
Recently, Rocchetta \emph{et al.} \cite{Rocchetta-Gentili-Sarfatti-MN-2012} established the following Bohr's theorem  for the class $\mathcal{SB}(\mathbb{B})$ of slice regular functions on $\mathbb{B}.$
\begin{theoB}{\rm \cite[Theorem 4.1, p.2101]{Rocchetta-Gentili-Sarfatti-MN-2012}}
	Let $f(q)=\sum_{k=0}^{\infty} q^kp_k\in \mathcal{SB}(\mathbb{B}).$ Then
	\begin{align}\label{Eq-1.111}
		\sum_{k=0}^{\infty}| q^kp_k|<1\;\;\mbox{for all }\;\; |q|\leq \dfrac{1}{3}.
	\end{align}
	Moreover, $1/3$ is the largest radius for which the statement is true.
\end{theoB}

In this paper, first our aim is to establish the Bohr inequality for the class $\mathcal{S}^*$ of slice starlike functions and the class $\mathcal{C}$ of slice close-to convex functions over quaternions. In addition, we establish results concerning Bohr phenomena for the class $\mathcal{SB}(\mathbb{B})$ of slice regular functions over quaternions.\vspace{2mm}

We organize the rest of the paper as follows: First we give the statements of our results for the classes $\mathcal{S}^*$, $\mathcal{C}$, and $\mathcal{SB}(\mathbb{B})$. In Section \ref{Sec-2}, we give details of preliminaries for these classes and quaternions. In Section \ref{Sec-3}, we give the details proof of the main results.\vspace{2mm}

We now state our first main result on Bohr’s theorem for the class $\mathcal{S}^*$.
\begin{thm}\label{Thm-1.1}
	Let $f(q)=q+\sum_{k=2}^{\infty}q^kp_k\in\mathcal{S}^*.$ Then
	\[\mathcal{A}_f(q):=|q|+\sum_{k=2}^{\infty}|q^kp_k|\leq 1\;\;\mbox{for}\;\; |q|=r\leq \frac{3-\sqrt{5}}{2}. \]
 The constant $(3-\sqrt{5})/2$ cannot be improved.
\end{thm}
The following result concerns Bohr's theorem for the class \( \mathcal{S}^* \), assuming the condition \( q f^{\prime}(q) \in \mathcal{S}^* \).
\begin{thm} \label{Thm-1.2}
	Let $f(q)=q+\sum_{k=2}^{\infty}q^kp_k$ be a slice regular function on $\mathbb{B}$ such that $qf^{\prime}(q)\in\mathcal{S}^*.$ Then
	\[\mathcal{B}_f(q):=|q|+\sum_{k=2}^{\infty}|q^kp_k|\leq 1\;\;\mbox{for}\;\; |q|\leq \dfrac{1}{2}. \]
	The constant $1/2$ cannot be improved.
\end{thm}
The following result we obtain the Bohr's theorem for the class $\mathcal{C}$.
\begin{thm}\label{Thm-1.3}
	Let $f(q)=q+\sum_{k=2}^{\infty}q^kp_k\in\mathcal{C}.$ Then
	\[|q|+\sum_{k=2}^{\infty}|q^kp_k|\leq 1\;\;\mbox{for}\;\; |q|=r\leq \frac{3-\sqrt{5}}{2}. \]
	The constant $(3-\sqrt{5})/2$ cannot be improved.
\end{thm}

In recent years, the generalization and refinement of Bohr-type inequalities have been an active area of research, with numerous studies focusing on their sharpness for specific classes of analytic functions. For a comprehensive overview, see \cite{Evdoridis-Ponnusamy-Rasila-RM-2021, Liu-Ponnusamy-PAMS-2021, Liu-Ponnusamy-Wang-RACSAM-2020, Liu-Liu-Ponnusamy-BSM-2021, Ponnusamy-Vijayakumar-Wirths-JMAA-2022, Allu-Arora-JMAA-2022} and references therein. Motivated by the results in \cite{Blasco-OTAA-2010}, we now pose the following question.

\begin{ques}\label{q-1.1}
	Can we extend the Bohr theorem to a more general form for the class $\mathcal{SB}(\mathbb{B})$?
\end{ques}
We define \(\mathcal{SB}_n(\mathbb{B})\) (\(n \geq 1\)) as the set of functions \(\omega \in \mathcal{SB}(\mathbb{B})\) satisfying
\begin{align*}
	\omega(0) = \omega'(0) = \dots = \omega^{(n-1)}(0) = 0\;\mbox{but}\; \omega^{(n)}(0)\neq 0.
\end{align*}
Using a generalized Schwarz lemma for slice regular functions over quaternions \cite[Theorem 3.10]{Xu-Ren-JGA-2018}, we obtain that   \begin{align}\label{Eq-1.4}
	|\omega(q)|\leq |q|^n\,\;\mbox{for}\;\mbox{all}\;\; q\in \mathbb{B}.
\end{align}
To address Question \ref{q-1.1} in the following result, we first derive a generalized form of Theorem B.
\begin{thm}\label{Thm-1.4}
	Let $m\in (0,2],$  $f(q)=\sum_{k=0}^{\infty} q^kp_k\in \mathcal{SB}(\mathbb{B})$ and $\omega_k \in \mathcal{SB}_k(\mathbb{B})$ $(k\geq 1).$ Then
	\begin{align*}
		\mathcal{K}_f(q):=|p_0|^m+\sum_{k=1}^{\infty}|\omega_k (q)p_k|\leq 1,\;\;\mbox{for}\;\; |q| \leq R_{m},
	\end{align*}
	where $R_{m}:={m}/{(2+m
		).}$ The constant $R_{m}$ is best possible for each $m$.
\end{thm}
\begin{rem}
	By setting \( m = 1 \) and \(\omega_k(q) = q^k\), Theorem \ref{Thm-1.4} simplifies to Theorem B, thereby confirming that Theorem \ref{Thm-1.4} generalizes the results of Theorem B.
\end{rem}
Motivated by \cite[Theorem 2]{Ponnusamy-Vijayakumar-Wirths-RM-2020}, we derive a sharper refined version of Theorem B for slice regular functions over quaternions.
\begin{thm}\label{Thm-1.5}
	Let $m\in (0,1],$  $f(q)=\sum_{k=0}^{\infty} q^kp_k\in \mathcal{SB}(\mathbb{B})$ and $\omega_k \in \mathcal{SB}_k(\mathbb{B})$ $(k\geq 1).$ Then
	\begin{align*}
		\mathcal{L}_f(q):=\mathcal{K}_f(q)+\left(\dfrac{1}{1+|p_0|}+\dfrac{|\omega_1(q)|}{1-|\omega_1(q)|}\right)\sum_{k=1}^{\infty}|\omega_k(q)p_k|^2\leq 1, \;\;\mbox{for}\;\; |q|\leq R_m=\dfrac{m}{m+2}.
	\end{align*}
	The constant ${m}/{(m+2)}$ is best possible for each $m$.
\end{thm}

Let $f$ be holomorphic in $\mathbb D$, and for $0<r<1$,  let $\mathbb D_r=\{z\in \mathbb C: |z|<r\}$.
Let $S_r:=S_r(f)$ denote the planar integral
\begin{align*}
	S_r=\int_{\mathbb D_r} |f^\prime(z)|^2 d A(z).
\end{align*}
If the function $f\in \mathcal{H}(\mathbb{D})$ has Taylor's series expansion $f(z)=\sum_{k=0}^{\infty}a_kz^k $, then we obtain  (see \cite{Kayumov-Ponnusamy-CRACAD-2018})
\begin{align*}
	S_r= \pi\sum_{k=1}^\infty k|a_k|^2 r^{2k}.
\end{align*}
We recall a remark from \cite{Kayumov-Ponnusamy-CRACAD-2018} regarding the quantity $S_r$ which provides a clear understanding of its significance.
\begin{rem}\cite{Kayumov-Ponnusamy-CRACAD-2018}
	Let us remark that if $f$ is a univalent then $S_r$ is the area of the image of the subdisk $\mathbb{D}_r$ under the mapping $f$. In the case of multivalent function, $S_r$ is greater than the area of the image of the subdisk $\mathbb{D}_r$. This fact could be shown by noting that
	\begin{align*}
		S_r=\int_{\mathbb D_r}|f^{\prime}(z)|^2dA(w)=\int_{\mathbb D_r}\nu_f(w)dA(w)\geq \int_{\mathbb D_r}dA(w)=\mbox{Area}(f(\mathbb{D}_r)),
	\end{align*}
	where $\mathbb{D}_r:=\{z\in\mathbb{C} : |z|<r\}$ and $\nu_f(w)=\sum_{f(z)=w} 1$ denotes the counting function of $f$.
\end{rem}
The quantity \( S_r \) plays a crucial role in the study of the improved Bohr inequality. Several results on this Bohr inequality have been established for the class \( \mathcal{B}(\mathbb{D}) \) (see \cite{Ismagilov-Kayumov-Ponnusamy-JMAA-2020, Kayumov-Ponnusamy-CRACAD-2018}) and for certain classes of sense-preserving harmonic mappings in the unit disk (see \cite{Evdoridis-Ponnusamy-Rasila-IM-2019}). The motivation for exploring the Bohr phenomenon using suitable combinations over quaternions arises from the discussions above and the observations presented in the following remark on the improved Bohr inequality.
\begin{rem}\label{Rem-2.1}
	Ismagilov \emph{et al.} \cite{Ismagilov-Kayumov-Ponnusamy-JMAA-2020} remarked that for any function $F : [0, \infty)\to [0, \infty)$ such that $F(t)>0$ for $t>0$, there exists an analytic function $f : \mathbb{D}\to\mathbb{D}$ given by $f(z)=\sum_{n=0}^{\infty}a_nz^n$, for which the inequality
	\begin{align*}
		\sum_{n=0}^{\infty}|a_n|r^n+\frac{16}{9}\left(\frac{S_r}{\pi}\right)+\lambda\left(\frac{S_r}{\pi}\right)^2+F(S_r)\leq 1\;\; \mbox{for}\;\; r\leq\frac{1}{3}
	\end{align*}
	is false, where $\lambda$ is given in \cite[Theorem 1]{Ismagilov-Kayumov-Ponnusamy-JMAA-2020}.
\end{rem}
Considering Remark \ref{Rem-2.1}, it is noteworthy to observe that augmenting a non-negative quantity with the Bohr inequality does not yield the desired inequality for the class $\mathcal{H}\left(\mathbb{D}\right)$. However, this observation motivates us to study the Bohr inequality further for slice regular functions over quaternions with a suitable setting.\vspace{1.2mm}

To advance the study of the Bohr inequality, it is natural to consider the following question: \textit{Is it possible to establish an improved Bohr inequality for slice regular functions over quaternions?} To address the question, we introduce a quantity analogous to \( S_r \). Henceforth, we define \( S_{q}^* \) for the slice regular function  $	f(q) = \sum_{k=0}^{\infty} q^k p_k \in \mathcal{SB}(\mathbb{B})$ as follows:
\begin{align}\label{Eq-1.2}
	S_{q}^*:=\sum_{k=1}^{\infty}k|q^kp_k|^2.
\end{align}
This observation motivates us to study the Bohr inequality further  for the class of functions $f\in \mathcal{SB}(\mathbb{B})$ with a suitable setting of $S_{q}^*$ and arise the following question.\vspace{1.2mm}
\begin{ques}\label{q-1.2}
	For the class $\mathcal{SB}(\mathbb{B}),$ can we establish an improved version of the Bohr inequality by incorporating the quantity  $S_{q}^*?$ If so, what conditions or constraints on $	S_{q}^*$ might lead to such an improvement?
\end{ques}
To give an answer to Question \ref{q-1.2}, in the next result, we first establish an improved version of the Bohr inequality for the class $\mathcal{SB}(\mathbb{B}).$
To serve our purpose, we consider a polynomial of degree $ N $ in $w$ as follows
\begin{align}\label{Eq-1.3}
	Q_N(w):=d_1w+d_2w^2+\dots+d_Nw^N, \;\; \mbox{where} \;\;
	d_i\geq 0, i=1,2,\dotsm N,
\end{align}
and obtain the following new  improved version of the Bohr inequality.
\begin{thm}\label{BS-thm-1.3}
	Let $m\in (0,1],$  $f(q)=\sum_{k=0}^{\infty} q^kp_k\in \mathcal{SB}(\mathbb{B})$ and $\omega_k \in \mathcal{SB}_k(\mathbb{B})$ $(k\geq 1).$ Then
	\begin{align*}
		\mathcal{M}_f(q):=|p_0|^m+\sum_{k=1}^{\infty}|\omega_k(q)p_k|+Q\left(	S^*_{q}\right)\leq 1\;\;\;\mbox{for}\;\;\;|q|\leq \dfrac{m}{2+m},
	\end{align*}
	where the coefficients of the polynomial $Q_N$ satisfy the condition
	\begin{align}\label{Eq-1.6}
		L(d_1,\dots, d_N):&=8d_1M_m^2+6c_2d_2M_m^4+\dots+2(2N-1)c_Nd_NM_m^{2N}
		\\ \nonumber
		&\leq m,
	\end{align}
	with $M_m:=m(2+m)/(4m+4)$ and
	\[
	c_k:=\max_{x\in [0,1]}\left(x(1+x)^2(1-x^2)^{2k-2}\right),
	\quad k=2,3,\cdots, N.
	\]
	The constant $ m/(2+m)$ is best possible for each $m$ and $d_1,\dots, d_N$ which satisfy $(\ref{Eq-1.6})$.
\end{thm}
Note that inequality \eqref{Eq-1.11} also can be written as
\begin{align*}
	\sum_{k=1}^{\infty}|a_k||z|^k\leq 1-|a_0|=1-|f(0)|=\mbox{dist} \left(f(0), \partial\mathbb{D}\right)\;\mbox{for all }\; |z|\leq \dfrac{1}{3},
\end{align*}
where `dist' is the Euclidean distance.\vspace{2mm}

In view of this distance formulation, many researchers studied the Bohr inequality in various settings for different classes of functions (see \cite{Liu-Ponnusamy-Wang-RACSAM-2020,Ponnusamy-Vijayakumar-Wirths-JMAA-2022,Xu-PRSE-2021}).  Inspired by these articles and there results, Xu and Ren \cite{Xu-Ren-JGA-2018} studied the Bohr inequality in a more general settings for slice regular functions over quaternions and obtained the following inequality.
\begin{theoC}{\rm \cite[Theorem 6.1, p.3800]{Xu-Ren-JGA-2018}}
	Let $ f(q)=\sum_{k=0}^{\infty}q^kp_k$ be slice regular function on $\mathbb{B}$ such that $f(q)\in \Pi:=\{ q\in \mathbb{H}: {\rm Re}(q)\leq 1\}$ for all $q\in \mathbb{B}.$ Then
\begin{align*}
	\sum_{k=1}^{\infty}|q^kp_k|\leq \mbox{dist} \left(f(0), \partial\Pi\right)\;\mbox{for}\; |q|\leq \dfrac{1}{3}.
\end{align*}
\end{theoC}

Now, our current objective is to address several additional questions related to the Bohr inequality. For instance, it is well known that the Bohr radius of \( {1}/{3} \) remains valid in the classical Bohr inequality even when the usual assumption on \( f \) is replaced by the condition \( \operatorname{Re} f(z) < 1 \) in \( \mathbb{D} \), with \( a_0 = f(0) \in [0,1) \). In fact, under this condition, it follows from Carathéodory's Lemma (see \cite[p.~41]{Duren-1983}) that \( |a_k| \leq 2(1 - a_0) \) for all \( k \in \mathbb{N} \). This leads to the following sharp inequality, as noted in \cite{Paulsen-Popescu-Singh-PLMS-2002}.
\begin{align*}
	|a_0|+\sum_{k=1}^{\infty}|a_k||z|^k\leq a_0+2(1-a_0)\dfrac{r}{1-r}\leq 1\;\;\mbox{for}\;\;|z|=r\leq \dfrac{1}{3}.
\end{align*}
Therefore, motivated by the article \cite[Theorem 1]{Ponnusamy-Vijayakumar-Wirths-JMAA-2022}, the following natural question arises for a slice regular function \( f \) in the open unit ball \( \mathbb{B} \) of \( \mathbb{H} \).

A natural question is to look for an analogue of \eqref{BS-Eq-1.2} if the image of $\mathbb{B}$ under the function $f$ in  $\Pi=\{ q\in \mathbb{H}: {\rm Re}(q)\leq 1\}.$ We answer this question in the next result.
\begin{thm}\label{Thm-1.7}
	Let $f$ be as in Theorem C and $p_0\in [0,1).$ Then
	\[ \mathcal{N}_f(q):=\sum_{k=0}^{\infty}|q^kp_k|+\left(\frac{1}{1+p_0}+\frac{|q|}{1-|q|}\right)\sum_{k=1}^{\infty}|q^{k}p_k|^2\leq 1\]
	holds for all $|q|=r\leq R_*$, where $R_*\approx 0.24683$ is the unique root in $(0, 1)$ of the equation $3r^3-5r^2-3r+1=0$. The constant $R_*$ is best possible.
\end{thm}

\section{\bf Preliminaries}\label{Sec-2}
In this section, we review essential definitions and fundamental results required for the subsequent discussion on slice regular functions.
\subsection{\bf The algebra of quaternions}
Let $\mathbb{H}$ denote the non-commutative, associative, real algebra of quaternions with standard basis $\{ 1, i,j,k\},$ subject to the multiplication rules
\begin{align*}
	i^2=j^2=k^2=ijk=-1.
\end{align*}
Every element $q=x_0+x_1i+x_2j+x_3k$ in  $\mathbb{H}$ is composed by the real part $\mathrm{Re}(q)=x_0$ and the imaginary part $\mathrm{Im}(q)=x_1i+x_2j+x_3k.$ The conjugate of $q\in \mathbb{H} $ is then $\overline{q}=\mathrm{Re}(q)-\mathrm{Im}(q)$ and its modulus is defined by $|q|^2=q\overline{q}=|\mathrm{Re}(q)|^2+|\mathrm{Im}(q)|^2.$ Therefore, we can calculate the multiplicative inverse of each $q\neq 0$ as $ q^-1=|q|^{-2}\overline{q}.$ Every  $q\in \mathbb{H} $ can be expressed as $q=x+yI,$ where $x,y\in \mathbb{R}$ and
\begin{align*}
	I=\dfrac{\mathrm{Im}(q)}{|\mathrm{Im}(q)|}
\end{align*}
if $\mathrm{Im}(q)\neq 0,$ otherwise we take $I$ arbitrarily such that $I^2=-1,$ where $I$ is an element of the unit $2$-sphere of purely imaginary quaternions
\[ \mathbb{S}:=\{q\in \mathbb{H}: q^2=-1 \}.\]
For every $I\in \mathbb{S},$ we will denote by $ \mathbb{C}_I$ the plane $\mathbb{R}\oplus I \mathbb{R},$ isomorphic to $\mathbb{C},$ and, if $\Omega\subset \mathbb{H},$ by $\Omega_I$ the intersection $\Omega\cap \mathbb{C}_I.$ \vspace{2mm}

We now revisit the definition of slice regularity.
\begin{defi}
	Let $\Omega$ be a domain in $\mathbb{H}$. A function $f :\Omega \to \mathbb{H}$ is called (left) slice regular if, for all $I \in \mathbb{S} $, its restriction $f_I$ to $\Omega_I$ is holomorphic, i.e., it has continuous partial derivatives and satisfies
	\[ \overline{\partial}_If(x+yI):=\dfrac{1}{2}\left(\dfrac{\partial}{\partial x}+I \dfrac{\partial}{\partial y}\right)f_I(x+yI)=0\]
	for all $x+yI\in \Omega_I.$
\end{defi}
\begin{defi}
Let $f:\mathbb{B} \to \mathbb{H}$ be a slice regular function. For each  $I \in \mathbb{S},$ the $I$-derivative of $f$ at $q=x+yI$ is defined by
\[ \partial_If(x+yI):=\dfrac{1}{2}\left(\dfrac{\partial}{\partial x}-I \dfrac{\partial}{\partial y}\right)f_I(x+yI)\]
on $\mathbb{B}_I.$ The slice derivative of $f$ is the function $f^\prime$ defined by $ \partial_If$ on $\mathbb{B}_I$ for all $I \in \mathbb{S}.$
\end{defi}
All  slice  regular function  on $\mathbb{B}$ can be expressed as power  series.

\begin{prop}\label{Th-1.1}
	Let $f:\mathbb{B} \to \mathbb{H}$ be a slice regular function. Then
	\[ f(q)=\sum_{k=0}^{\infty} q^kp_k,\;\;\;\mbox{with}\;\; p_k=\frac{f^{k}(0)}{k!}\] for all $q\in \mathbb{B}.$
\end{prop}

In fact, there are some limitations in the expansion given in Proposition \ref{Th-1.1}. Hence, a new type of expansion occurs in \cite{Stoppato-AM-2012}.

\begin{defi}
	Let $f,h:\mathbb{B} \to \mathbb{H}$ be two slice regular functions of the form
	\[ f(q)=\sum_{k=0}^{\infty} q^kp_k,\;\;\; h(q)=\sum_{k=0}^{\infty} q^ka_k.\]
	 The  regular product (or $\ast$-product) of $f$ and $h$ is the slice regular function defined by
	 \[ f\ast h(q)=\sum_{k=0}^{\infty } q^k\left( \sum_{n=0}^{k}p_na_{k-n}\right).\]
\end{defi}
The relationship between the $\ast$-product and the usual pointwise product is explained by the following result.
\begin{prop}
		Let $f,h:\mathbb{B} \to \mathbb{H}$ be two slice regular functions. Then for all $q \in \mathbb{B},$
		\[ f\ast h(q)=	\begin{cases}
			f(q)h(f(q)^{-1}qf(q)), \;\;\;\mbox{if}\;\; f(q)\neq 0, \vspace{2mm}\\
			0, \;\;\;\;\;\;\;\;\;\;\;\;\;\;\;\;\;\;\;\;\;\;\;\;\;\;\;\;\;\;\;\;\mbox{if}\;\; f(q)= 0.
		\end{cases}\]
\end{prop}
\begin{defi}
Let $f(q)=\sum_{k=0}^{\infty} q^kp_k$ be a slice regular function on $\mathbb{B}.$  We define the regular conjugate of $f$ as
\[f^c(q)=\sum_{k=0}^{\infty} q^k\overline{p}_k, \]
and the symmetrization of $f$ as
\[ f^s(q)=f\ast f^c(q)=f^c\ast f(q)=\sum_{k=0}^{\infty } q^k\left( \sum_{n=0}^{k}p_n\overline{p}_{k-n}\right). \]
\end{defi}
Both $f^c$ and $f^s$ are slice regular functions on $\mathbb{B}.$\vspace{1.2mm}

Let $Z_{f^s}$ denote the zero set of the symmetrization $f^s$ of $f.$
\begin{defi}
	Let $f$ be a slice regular function on $\mathbb{B}.$ If  $f$ does not vanish  identically, its regular reciprocal is the function defined by
	\[ f^{-\ast}(q):= f^s(q)^{-1}f^c(q)\]
	 which is slice regular on $\mathbb{B}\setminus Z_{f^s}.$
\end{defi}
\begin{prop}
	Let $f$ and $g$ be slice regular functions on $\mathbb{B}.$ Then for all  $q \in \mathbb{B}\setminus Z_{f^s},$
	\[  f^{-\ast}(q) \ast g(q)= f\left(T_f(q)\right)^{-1}g\left( T_f(q)\right),\]
	 where $T_f:\mathbb{B}\setminus Z_{f^s} \to \mathbb{B}\setminus Z_{f^s}$ is defined by $T_f(q)=f^c(q)^{-1}qf^c(q).$ Furthermore, $T_f$ and $T_{f^c}$ are mutual inverses so that $T_f$ is a diffeomorphism.
\end{prop}
We state the definition of slice starlike  functions. Let $f$ be a slice regular function in the unit ball $\mathbb{B}$ of quaternions $\mathbb{H},$ normalized by $f(0)=0$ and $f^\prime(0)=1.$ We further assume that $Z_f=\{ 0\},$ where $Z_f$ denotes the zero set of $f.$ Denote this function class by $\mathcal{S}.$  Notice that this condition is weaker than $f$ is injective. The set of slice starlike functions and slice  close-to-convex functions are defined respectively as follows:
\[ \mathcal{S}^\ast:=\{f:\mathbb{B} \to \mathbb{H}: f\in \mathcal{S}\; \mbox{with}\; {\rm Re}\left(f(q)^{-1}qf^{\prime}(q)\right)>0\;\mbox{on}\;\mathbb{B}  \}\] and
\[ \mathcal{C}=\{ f:\mathbb{B} \to \mathbb{H}: \exists\; h\in \mathcal{S}^\ast\; \mbox{such \; that }\; {\rm Re}\left(h(q)^{-1}qf^{\prime}(q)\right)>0\;\mbox{on}\;\mathbb{B} \}.\]

 \section{\bf Proof of main results}\label{Sec-3}
 \begin{proof}[\bf Proof of Theorem \ref{Thm-1.1}]
 	Let $f(q)=q+\sum_{k=2}^{\infty}q^kp_k\in\mathcal{S}^*.$ It follows from \cite[Theorem 1.1]{Xu-Ren-JGA-2018} that
 	\begin{align}\label{Eq-1.8}
 		|p_k|\leq k\;\;\;\mbox{for}\;\;\; k=2,3,\dots.
 	\end{align}
 	In view of \eqref{Eq-1.8}, by a straightforward calculation, we obtain
 	\begin{align*}
 		\mathcal{A}_f(q)&\leq  |q|+\sum_{k=2}^{\infty}k|q|^k= |q|+\dfrac{|q|^2(2-|q|)}{(1-|q|)^2}\\&=1+\dfrac{|q|(1-|q|)^2+(2-|q|)|q|^2-(1-|q|)^2}{(1-|q|)^2}\\&=1+\dfrac{(-1+3|q|-|q|^2)}{(1-|q|)^2}\leq 1\;\;\mbox{for}\;\; |q|\leq r_*,
 	\end{align*}
 	where $r_*$ is the unique root in $(0,1)$ of the equation $-1+3r-r^2=0.$ Thus the desired inequality is established.\vspace{2mm}
 	
 	To prove the sharpness part of the theorem, for some $u\in \partial \mathbb{B},$ we consider the function $f\in \mathcal{S}^* $ given by
 	\begin{align}
 		f_u(q)=q(1-qu)^{-*2}=\sum_{k=1}^{\infty} q^kp_k, \;\; q\in \mathbb{B},
 	\end{align}
 	where	$p_k=ku^{k-1},$ $k\in \mathbb{N}.$ For the function $f_u,$ a simple computation leads to
 	\begin{align*}
 		\mathcal{A}_{f_u}(q)&=|q|+\sum_{k=2}^{\infty}|kq^ku^{k-1}|=|q|+\frac{(2-|q|)|q|^2}{(1-|q|)^2}\\&=1+\frac{-1+3|q|-|q|^2}{(1-|q|)^2}
 =1+\frac{(|q|-r_*)((3+\sqrt{5})/2-|q|)}{(1-|q|)^2}.
 	\end{align*}
 	It is evident that $\mathcal{A}_{f_u}(q) > 1$ for $|q| > r_*$, 
 which confirms that the constant $r_*$ cannot be improved.
 \end{proof}

 \begin{proof}[\bf Proof of Theorem \ref{Thm-1.2}]
 	Under the condition of Theorem \ref{Thm-1.2}, it follows that (see \cite[Theorem 1.3]{Xu-Ren-JGA-2018})
 	\begin{align}\label{Eq-1.10}
 		|p_k|\leq 1\;\;\mbox{for}\;\; k=2,3,\dots.
 	\end{align}
 	In view of \eqref{Eq-1.10}, by a straightforward calculation, we obtain
 	\begin{align*}
 		\mathcal{B}_f(q)&\leq |q|+\sum_{k=2}^{\infty}|q|^k\\&=|q|+\dfrac{|q|^2}{1-|q|}=1+\dfrac{2|q|-1}{1-|q|}.
 	\end{align*}
 	It is easy say that $\mathcal{B}_f(q)\leq 1$ for $|q|\leq 1/2.$\vspace{1.2mm}
 	
 	To prove the sharpness, for some $u\in \partial\mathbb{B}$, we consider the function
 	\[ f_u(q)=q(1-qu)^{-*}=q+\sum_{k=2}^{\infty}q^ku^{k-1},\;\; q\in \mathbb{B}.\]
 	For the function $f_u,$ we obtain
 	\begin{align*}
 		\mathcal{B}_{f_u}(q)&=|q|+\sum_{k=2}^{\infty}|q^ku^{k-1}|\\&=|q|+\dfrac{|q|^2}{1-|q|}=1+\dfrac{2|q|-1}{1-|q|}.
 	\end{align*}
 	This shows that $\mathcal{B}_{f_u}(qu)>1$ for $|q|>1/2.$ This completes the proof.
 \end{proof}

 \begin{proof}[\bf Proof of Theorem \ref{Thm-1.3}]
 	Under the condition of Theorem \ref{Thm-1.2}, it follows that (see \cite[Theorem 4.3]{Xu-Ren-JGA-2018})
 	\begin{align*}
 		|p_k|\leq k\;\;\mbox{for}\;\; k=2,3,\dots.
 	\end{align*}
 	The rest of the proof follows from Theorem \ref{Thm-1.1} and therefore, we omit the details.
 \end{proof}

\begin{proof}[\bf Proof of Theorem \ref{Thm-1.4}]
		Let $f(q)=\sum_{k=0}^{\infty} q^kp_k\in \mathcal{SB}(\mathbb{B}).$  	Under the condition of Theorem \ref{Thm-1.4}, it follows that (see \cite[Theorem 4.1]{Rocchetta-Gentili-Sarfatti-MN-2012})
	\begin{align}\label{Eq-2.4}
		|p_k|\leq 1-|p_0|^2,\;\;\mbox{for}\;\;k\in \mathbb{N}.
	\end{align}
	In view of  \eqref{Eq-1.4} and \eqref{Eq-2.4}, we obtain
	\begin{align*}
		\mathcal{K}_f(q)&\leq|p_0|^m+(1-|p_0|^2)\dfrac{|q|}{\left(1-|q|\right)}\leq 1
	\end{align*}
	if, and only if, $|q|\leq R_{m},$ where
	\begin{align*}
		R_{m}= \inf_{|p_0|<1} \{\mathcal{F}\left(|p_0|\right)\}\;\;\mbox{and}\;\;   \mathcal{F}\left(|p_0|\right):=\dfrac{1-|p_0|^m}{2-|p_0|^2-|p_0|^m}.
	\end{align*}
	Moreover, for $0<m\leq 2,$ a computations gives that
	\begin{align*}
		\mathcal{F}^{\prime}(t)=\dfrac{t^{m-1}B(t)}{\left(2-t^2-t^m\right)^2},\;\;\mbox{where}\;\;B(t):=-m-(2-m)t^2+2t^{2-m}.
	\end{align*}
	Since $B^\prime(t)=2(2-m)t^{1-m}(1-t^m)\geq 0$ for $m\in (0,2]$ and $0< t\leq 1$, this implies that $B(t)\leq B(1)=0$ and, hence $\mathcal{F}$ is a decreasing function of $t\in [0,1).$ Therefore,
	\begin{align*}
		\mathcal{F}(t)\geq \lim\limits_{t\rightarrow 1^{-}}\mathcal{F}(t)=\dfrac{m}{2+m}=R_{m}\, \mbox{ for } t\in [0,1).
	\end{align*}
	This proves the desired inequality.\vspace{1.2mm}
	
	Finally, to show the sharpness, given $a\in (0,1)$ and $u\in \partial \mathbb{B},$ we consider the slice regular function
	\begin{align}\label{Eq-1.1}
		f_a(q)=(1-qa)^{-\ast }\ast(a-q)u=a-(1-a^2)u\sum_{k=1}^{\infty}q^ka^{k-1},\;\;q\in \mathbb{B}.
	\end{align}
	For $\omega_k(q)=q^ku$ and $f_a,$ we obtain
	\begin{align*}
		\mathcal{K}_{f_a}(q)&=a^m+(1-a^2)\sum_{k=1}^{\infty}a^{k-1}|q|^k\\&=a^m+\dfrac{(1-a^2)|q|}{1-a|q|}\\&=1+(1-a)\mathcal{G}(a),
	\end{align*}
	where
	\begin{align*}
		\mathcal{G}(a):=\dfrac{(1+a)|q|}{1-a|q|}-\left(\dfrac{1-a^m}{1-a}\right).
	\end{align*}
	Taking $a$ very closed to $1$, that is
	\begin{align*}
		\lim\limits_{a\rightarrow 1^{-}}\mathcal{G}(a)=\dfrac{2 |q|}{1-|q|}-m=\dfrac{(2+m) |q|-m}{1-|q|}>0\, \mbox{ for }\, |q|>R_{m}=\dfrac{m}{2+m}.
	\end{align*}
	Therefore, $\mathcal{K}_{f_a}(q)>1$ for  $|q|>R_m$, which shows that the constant $R_m$ is best possible and the proof is completes.
\end{proof}

 \begin{proof}[\bf Proof of Theorem \ref{Thm-1.5}]
 	Let $f(q)=\sum_{k=0}^{\infty} q^kp_k\in \mathcal{SB}(\mathbb{B}).$ In view of \eqref{Eq-2.4}, we obtain
 	\begin{align*}
 		\mathcal{L}_f(q)&\leq |p_0|^m+(1-|p_0|^2)\sum_{k=1}^{\infty}|q|^k+\left(\dfrac{1}{1+|p_0|}+\dfrac{|q|}{1-|q|}\right)(1-|p_0|^2)^2\sum_{k=1}^{\infty}|q|^{2k}\\&=1+(|p_0|^m-1)+
 (1-|p_0|^2)\dfrac{|q|}{1-|q|}+\left(\dfrac{1}{1+|p_0|}+\dfrac{|q|}{1-|q|}\right)\dfrac{(1-|p_0|^2)^2|q|^2}{1-|q|^2}\\&:=\mathcal{L}^*_f(q).
 	\end{align*}
 	To proceed further in the proof, we use the following inequality proved in \cite{Lin-Liu-Ponnu-ACS-2021}
 	\begin{align}\label{Eq-2.5}
 		 \frac{1-t^m}{1-t}\geq m\; \mbox{for}\; m\in (0, 1]\; \mbox{and}\; t\in [0, 1).
 	\end{align}
 	In view of \eqref{Eq-2.5}, we obtain that
 	\begin{align*}
 		\mathcal{L}^*_f(q)\leq 1+m(t-1)+K_1(1-t^2)+K_2(1-t)(1-t^2)+K_3(1-t^2)^2:=\Psi(t),
 	\end{align*}
 	where $t=|p_0|\in [0,1],$ and
 	\begin{align*}
 		K_1:=\dfrac{|q|}{1-|q|},\;\;K_2:=\dfrac{|q|^2}{1-|q|^2}\;\;\mbox{and}\;\; K_3:=\dfrac{|q|^3}{(1-|q|)(1-|q|^2)}.
 	\end{align*}
 	Clearly, $K_1,$ $K_2$ and $K_3$ are non-negative for all $q\in \mathbb{B}\setminus\{0\}$. Also, We see that
 	\begin{align*}
 		\begin{cases}
 			\Psi^\prime(t)=m-2K_1t+K_2(3t^2-2t-1)+4K_3(t^3-t),
 			\vspace{2mm}\\
 			\Psi^{\prime\prime}(t)=-2K_1+2K_2(3t-1)+4K_3(3t^2-1)
 			\vspace{2mm}\\
 			\Psi^{\prime\prime\prime}(t)=6K_2+24tK_3.
 		\end{cases}
 	\end{align*}
 	Since $K_2$ and $K_3$ are non-negative, we have $\Psi^{\prime\prime\prime}(t)>0$ for all $t\in [0,1].$ Thus, $	\Psi^{\prime\prime}$ is an increasing function of $t,$ which implies that
 	\begin{align*}
 		\Psi^{\prime\prime}(t)\leq \Psi^{\prime\prime}(1)&=-2K_1+4K_2+8K_3\\&= \dfrac{2|q|\left( 4|q|^2+2|q|(1-|q|)-(1-|q|^2) \right)}{(1-|q|)(1-|q|^2)}\\&=\dfrac{2|q|(1+|q|)(3|q|-1)}{(1-|q|)(1-|q|^2)}.
 	\end{align*}
 	It is easy to see that $\Psi^{\prime\prime}(t)\leq0$ for $|q|\leq 1/3$ and $t\in [0,1],$ which shows that $\Psi^{\prime}$ is a decreasing function in $[0,1].$ Therefore, we obtain
 	\begin{align*}
 		\Psi^{\prime}(t)\geq \Psi^{\prime}(1)=m-2K_1=\dfrac{m-(2+m)|q|}{1-|q|}\geq 0\;\;\mbox{for}\;\; |q|\leq \dfrac{m}{2+m}.
 	\end{align*}
 	Hence, $\Psi^{\prime}(t)\geq 0$ for $t\in [0,1]$ and $|q|\in \left(0,m/(m+2)\right),$ which shows that $\Psi^{\prime}$ is an increasing function in $[0,1].$ Thus, we obtain
 	\begin{align*}
 		\Psi(t)\leq \Psi(1)=1\;\;\mbox{for}\;\; |q|\leq \dfrac{m}{2+m}.
 	\end{align*}
 	This proves the desired inequality.\vspace{1.2mm}
 	
 	To show the constant ${m}/{(2+m)}$ is best possible, given $a\in (0,1),$ we consider slice regular function $f_a$ given by \eqref{Eq-1.1}. For $\omega_k(q)=q^ku$ (k$\geq$ 1) for some $u\in \partial\mathbb{B}$ and $f_a,$ we obtain
 	\begin{align*}
 		\mathcal{L}_{f_a}(q)&=a^m+(1-a^2)\sum_{k=1}^{\infty}|q|^ka^{k-1}+\left(\dfrac{1}{1+a}+\dfrac{|q|}{1-|q|}\right)(1-a^2)^2\sum_{k=1}^{\infty}|q|^{2k}a^{2(k-1)}\\&= a^m+\dfrac{(1-a^2)|q|}{ 1-|q|}=1+\frac{(2-a^2-a^m)|q|-(1-a^m)}{1-|q|}>1\;\;\mbox{for}\;\;|q|> \dfrac{m}{2+m}.
 	\end{align*}
 	This implies that the constant ${m}/{(2+m)}$ is best possible. This completes the proof.
 \end{proof}

 \begin{proof}[\bf Proof of Theorem \ref{BS-thm-1.3}]
 	Let $f(q)=\sum_{k=0}^{\infty} q^kp_k\in \mathcal{SB}(\mathbb{B}).$ In view of \eqref{Eq-1.2} and \eqref{Eq-2.4}, we obtain
 	\begin{align}\label{Eq-1.5}
 		S^*_{q}\leq (1-|p_0|^2)^2\sum_{k=1}^{\infty}k |q|^{2k}=(1-|p_0|^2)^2\dfrac{|q|^2}{\left(1-|q|^2\right)^2}\;\;\mbox{for}\;\;|q|<1.
 	\end{align}
 	
 	Let $|p_0|=t.$  In view of \eqref{Eq-1.3}, \eqref{Eq-1.5}, \eqref{Eq-2.4}, and \eqref{Eq-2.5},  by a simple computation, the following inequality can be obtained
 	\begin{align*}
 		\mathcal{M}_f(q)&\leq 1+ m(t-1)+ \left(1-t^2\right)\frac{|q|}{1-|q|}+\sum_{k=1}^{N}d_k\left(\dfrac{(1-t^2)|q|}{1-|q|^2}\right)^{2k}= 1+\Phi(t,|q|),
 	\end{align*}
 	where \begin{align*}
 		\Phi(t,|q|):=\dfrac{(1-t^2)|q|}{1-|q|}+\sum_{k=1}^{N}d_k\left(\dfrac{(1-t^2)|q|}{1-|q|^2}\right)^{2k}-m(1-t).
 	\end{align*}
 	For all $t\in [0,1)$, by a straightforward calculation, it can be shown that $\Phi(t,|q|)$ is a monotonically increasing function of $|q|$. Consequently, for $t\in [0,1)$ $|q|\leq \dfrac{m}{2+m}$, we obtain
 	\[ \Phi(t,|q|)\leq \Phi(t,m/(2+m)).\]
 	 By a straightforward calculation
 	\begin{align*}
 		\Phi(t,m/(2+m))=\dfrac{(1-t^2)}{2}\left(m+2F_N(t)-\dfrac{2m}{1+t}\right)=\dfrac{(1-t^2)}{2}\mathcal{W}(t),
 	\end{align*}
 	where $M_m=\frac{m((m+2)}{4(m+1)}$, and
 \begin{align*}
 		F_N(t):=\sum_{k=1}^{N}d_k(1-t^2)^{2k-1}M^{2k}_m\;\;\mbox{and}\;\; \mathcal{W}(t):=m+2F_N(t)-\dfrac{2m}{1+t}.
 	\end{align*}
 	To establish $\Phi(t,|q|)\leq 0,$ it suffices to show that $\mathcal{W}(t)\leq 0$ for $t\in [0,1].$ As $t\in [0,1 ]$, a straightforward calculation reveals that
 	\begin{align*}
 		t(1+t)^2M_m^2&\leq 4 M_m^2,\\ t(1+t)^2(1-t^2)^2M_m^4&\leq c_2M_m^4,\\ \vdots \\t(1+t)^2(1-t^2)^{2N-2}M_m^{2N}&\leq c_NM_m^{2N}.
 	\end{align*}
 	Thus, we see that
 	\begin{align*}
 		\mathcal{W}^{\prime}(t)&=\frac{2}{(1+t)^2}\bigg(m-2d_1t(1+t)^2M_m^{2}-6d_2t(1+t)^2(1-t^2)^2M_m^{4}-\cdots\\&\quad-2(2N-1)d_Nt(1+t)^2(1-t^2)^{2N-2}M_m^{2N}\\&\geq \dfrac{2}{(1+t)^2} (m-L(d_1,\cdots,d_N))\\&\geq 0,
 	\end{align*}
 	if  the coefficients $d_i$ of the polynomial $Q_N$ satisfy the condition outlined in $(\ref{Eq-1.6})$. This indicates that $\mathcal{W}(t)$ behaves as an ascending function in $t\in [0,1]$, leading to the conclusion that $\mathcal{W}(t)\leq\mathcal{W}(1)=0$. This, in turn, establishes the desired inequality.\vspace{1.2mm}
 	
 	To show the sharpness, given $a\in (0,1),$ we consider slice regular function $f_a$ given by \eqref{Eq-1.1} and $\omega_k(q)=q^ku$ for some $u\in \partial\mathbb{B}$. In view of \eqref{Eq-1.2} and \eqref{Eq-1.3}, we obtain
 	\begin{align*}
 		\mathcal{M}_{f_a}(q)=1-(1-a)\Psi^*(|q|,a),
 	\end{align*}
 	where
 	\begin{align*}
 		\Psi^*(|q|,a)&:=\dfrac{1-|q|^m}{1-|q|}-\dfrac{(1+a)|q|}{1-|q|}-\dfrac{d_1|q|^2(1-a)(1+a)^2}{(1-a^2|q|^2)^2}-\dots\\&\quad-\dfrac{d_N|q|^{2N}(1-a)^{2N-1}(1+a)^{2N}}{(1-a^2|q|^2)^{2N}}.
 	\end{align*}
 	Therefore, we have
 	\[
 	\lim\limits_{a\to 1^{-}}\Psi^*(|q|,a)=m-\frac{2|q|}{1-|q|}<0\;\mbox{for}\;\;|q|>\dfrac{m}{2+m}.
 	\]
 	Thus, it follows that $\Psi^*(|q|,a)<0$ for $a$ sufficiently close to $1$.
 	Hence, we have
 	\begin{align*}
 		\mathcal{M}_{f_a}(q)=1-(1-a)\Psi^*(|q|)>1,
 	\end{align*}
 	which shows that ${m}/{(2+m)}$ is best possible and thereby concluding the proof.
 \end{proof}

 \begin{proof}[\bf Proof of Theorem \ref{Thm-1.7}]
 	Let $\alpha=1-p_0$. Under the condition of Theorem \ref{Thm-1.7}, it follows that (see \cite[p.1231]{Ren-Wang-CAOT-2015})
 	\begin{align}\label{Lem-2.11}
 		|p_k|\leq 2\left(1-{\rm Re}(f(0))\right),\;\;\;\mbox{for}\;\;\; k\in\mathbb{N}.
 	\end{align}
 	In view of \eqref{Lem-2.11}, we obtain $|p_k|\leq 2\alpha$ for all $k\geq 1$. For $\alpha\in (0, 1]$ and $|q|=r<1,$ we obtain
 	\begin{align*}
 		\mathcal{N}_f(q)&\leq 1-\alpha+\frac{2\alpha r}{1-r}+\left(\frac{1}{1+p_0}+\frac{r}{1-r}\right)\frac{4\alpha^2r^2}{1-r^2}\\&=1-\alpha\left(\frac{1-3r}{1-r}-\left(\frac{1+r-\alpha r}{(1-r)(2-\alpha)}\right)\frac{4\alpha r^2}{1-r^2}\right)\\&=1-\alpha\frac{Q(\alpha, r)}{(2-\alpha)(1-r)(1-r^2)},
 	\end{align*}
 	where \[ Q(\alpha, r):=4r^3\alpha^2-(7r^3+3r^2-3r+1)\alpha+6r^3-2r^2-6r+2.\]
 Our aim is to show that $ Q(\alpha, r)\geq 0$ for $r\leq R_*$ and $\alpha\in (0, 1]$. It follows that $\frac{\partial^2 Q(\alpha, r)}{\partial \alpha^2}\geq 0$ for any $\alpha\in (0, 1]$ and thus, $\frac{\partial Q}{\partial \alpha}$ is an increasing function of $\alpha$. This gives
 	\begin{align*}
 		\frac{\partial Q(\alpha, r)}{\partial \alpha}\leq \frac{\partial Q(1, r)}{\partial\alpha}=-(1-r)^3<0\;\;\mbox{for}\;\;r<1.
 	\end{align*}
This shows that $Q$ is a decreasing function of $\alpha$ on $(0, 1]$ so that
 	\begin{align*}
 		Q(\alpha, r)\geq Q(1, r)=3r^3-5r^2-3r+1,
 	\end{align*}
 	which is greater than or equal to $0$ for all $r\leq R_*$, where $R_*\approx 0.24683$ is the unique root of the equation $3r^3-5r^2-3r+1=0$ which lies in $(0, 1)$. \vspace{1.2mm}
 	
 	To prove the sharpness, for $u\in \partial \mathbb{B}$ and $a\in (0,1),$ we consider the function
 	\[ g_a(q)=a-2(1-a)q(1-q)^{-\ast}u=a-\sum_{k=1}^{\infty}q^k2(1-a)u.\]
 	For $g_a,$ we obtain
 	\begin{align*}
 		 \mathcal{N}_{g_a}(q)&=\sum_{k=0}^{\infty}|q^kp_k|+\left(\frac{1}{1+p_0}+\frac{|q|}{1-|q|}\right)\sum_{k=1}^{\infty}|q^{k}p_k|^2\\&=a+2(1-a)\dfrac{|q|}{1-|q|}+\left(\frac{1}{1+a}+\frac{|q|}{1-|q|}\right)\dfrac{4(1-a)^2|q|^2}{1-|q|^2}.
 	\end{align*}
 By the similar computations, it can be easily shown that $\mathcal{N}_{g_a}(q)>1$ for $|q|>R_*.$ This completes the proof of the theorem.
 \end{proof}\vspace{2mm}

  \noindent{\bf Acknowledgment:} 
  The research of the first author is supported by UGC-JRF (Ref. No. 201610135853), New Delhi, Govt. of India and second author is supported by SERB File No. SUR/2022/002244, Govt. of India.\vspace{1.5mm}

 \noindent {\bf Funding:} Not Applicable.\vspace{1.5mm}

 \noindent\textbf{Conflict of interest:} The authors declare that there is no conflict  of interest regarding the publication of this paper.\vspace{1.5mm}

 \noindent\textbf{Data availability statement:}  Data sharing not applicable to this article as no datasets were generated or analysed during the current study.\vspace{1.5mm}

 \noindent{\bf Code availability:} Not Applicable.\vspace{1.5mm}

 \noindent {\bf Authors' contributions:} All the authors have equal contributions to prepare the manuscript.

\end{document}